\documentclass[12pt]{amsart}

\newcommand{\Y}{{\mathcal Y}}
\newcommand{\LL}{{\mathcal L}}
\newcommand{\MM}{{\mathcal M}}

\newcommand{\Abar}{{\overline{A}}}
\newcommand{\fbar}{{\overline{f}}}
\newcommand{\Xbar}{{\overline{X}}}
\newcommand{\pibar}{{\overline{\pi}}}

\newcommand{\C}{{\mathbf C}}
\newcommand{\Q}{{\mathbf Q}}

\newcommand{\kbar}{{\overline{k}}}

\newcommand{\Z}{{\mathbf Z}}
\newcommand{\R}{{\mathbf R}}
\newcommand{\PP}{{\mathbf P}}

\newcommand{\G}{{\mathbf G}}
\newcommand{\tors}{_{\operatorname{tors}}}

\newcommand{\Gal}{\operatorname{Gal}}

\newcommand{\Pic}{\operatorname{Pic}}

\newcommand{\isom}{\cong}
\newcommand{\tensor}{\otimes}

\newtheorem{theorem}{Theorem$\!\!$}	
\newtheorem{lemma}[theorem]{Lemma$\!\!$}	

\newtheorem{prop}[theorem]{Proposition$\!\!$}	
\newtheorem{conj}{Conjecture$\!\!$}	

\theoremstyle{definition}

\theoremstyle{remark}
\newtheorem{rem}{Remark$\!\!$}	
\newtheorem{rems}{Remarks$\!\!$}	

\begin{document}

\title{Mordell-Lang plus Bogomolov}
\author{Bjorn Poonen}
\thanks{This research was supported by
a fellowship from the Gabriella and Paul Rosenbaum Foundation,
and by grants from the U.~C.~Berkeley Committee on Research.}
\address{Department of Mathematics, University of California,
		Berkeley, CA 94720-3840, USA}
\email{poonen@math.berkeley.edu}
\date{June 29, 1998}

\maketitle

\section{Theorem}

Let $k$ be a number field.
Let $A$ be an {\em almost split} semiabelian variety over $k$;
by this we mean that $A$ is isogenous to the product of
an abelian variety $A_0$ and a torus $T$.
We enlarge $k$ if necessary to assume that $T \isom \G_m^n$.
Let $\phi=(\phi_1,\phi_2):A \rightarrow A_0 \times \G_m^n$ be the isogeny.
Let $h_1: A_0(\kbar) \rightarrow \R$ be a N\'eron-Tate canonical height
associated to a symmetric ample line bundle on $A_0$,
and let $h_2: \G_m^n(\kbar) \rightarrow \R$ be the sum of the naive
heights of the coordinates.
For $x \in A(\kbar)$, let $h(x)=h_1(\phi_1(x)) + h_2(\phi_2(x))$.
For $\epsilon \ge 0$,
let $B_\epsilon=\{\, z \in A(\kbar) \mid h(z) \le \epsilon \,\}$.
Let $\Gamma$ be a finitely generated subgroup of $A(\kbar)$,
and define
	$$\Gamma_\epsilon := \Gamma+B_\epsilon =
	\{\, \gamma+z \mid \gamma \in \Gamma, h(z) \le \epsilon \,\}.$$
Note that $\Gamma_0 = \Gamma + A(\kbar)\tors$.

Let $X$ be a geometrically integral closed subvariety of $A$.
Our main result is the existence of $\epsilon>0$ such that
$X(\kbar) \cap \Gamma_\epsilon$ is contained
in a finite union $\bigcup Z_j$ where each $Z_j$ is a translate
of a sub-semiabelian variety of $A_\kbar=A \tensor_k \kbar$
by a point in $\Gamma_0$ and $Z_j \subseteq X_\kbar$.
Since a sub-semiabelian variety of an almost split
semiabelian variety is almost split,
this is equivalent to

\begin{theorem}
\label{main}
If $X$ is not a translate of a sub-semiabelian variety of $A$
by an element of $\Gamma_0$,
then there exists $\epsilon>0$ such that
$X(\kbar) \cap \Gamma_\epsilon$ is not Zariski dense in $X$.
\end{theorem}

Let us indicate the relation to more familiar statements.
The ``Mordell-Lang conjecture'' states that if $X$
is a geometrically integral closed subvariety of a semiabelian variety $A$
over a field $k$ of characteristic~0,
if $\Gamma$ is a finite rank subgroup of $A(\kbar)$,
and if $X$ is not the translate of a sub-semiabelian variety,
then $X(\kbar) \cap \Gamma$ is not Zariski dense in $X$.
The statement in this form was conjectured by Lang~\cite{lang}
and proved by McQuillan~\cite{mcquillan}
following work by Faltings, Vojta, Raynaud, Hindry, and many others.
See~\cite{mcquillan} and the references cited there for further history.
(There is a function field version too~\cite{hrushovski},
but it will not be considered here.)
The case where $A$ is almost split and where $\Gamma$
is finitely generated follows from our theorem,
since specialization arguments reduce the
Mordell-Lang conjecture to the case where $k$ is a number field.

The ``generalized Bogomolov conjecture''
states that if $A$ is an abelian variety over a number field $k$,
if $h$ is a N\'eron-Tate height,
and if $X$ is a geometrically integral closed subvariety
that is not the translate of a sub-abelian variety by a torsion point,
then $X(\kbar) \cap B_\epsilon$ is not Zariski dense in $X$.
This was proved by Zhang~\cite{zhangbogomolov},
following work by him and Szpiro and Ullmo,
and an independent proof was given shortly thereafter
by David and Philippon~\cite{davidphilippon}.
See~\cite{ullmo} and the survey by Abb\`es~\cite{abbes} for some history.
The analogue for tori was also proved by Zhang~\cite{zhangtori},
and a version for almost split semiabelian varieties
has been announced recently by Chambert-Loir~\cite{chambertloir2}.
These results are also contained in our theorem.

The proof of our theorem does not yield new proofs of any of these results,
however, because it uses them!
Our proof requires also the theorems on
``equidistribution of small points''
developed by Szpiro, Ullmo, Zhang~\cite{suz},\cite{zhangbogomolov},
Bilu~\cite{bilu}, and Chambert-Loir~\cite{chambertloir2}.

In the final section of our paper, we formulate a conjecture
which would generalize our theorem to all semiabelian varieties,
and with the division group of $\Gamma$ in place of $\Gamma$.

\section{Proof}

\begin{lemma}
\label{measures}
Let $V$ be a projective variety over $\C$.
Let $S$ be a connected quasi-projective variety over $\C$.
Let $\Y \rightarrow V \times S$ be a closed immersion of $S$-varieties,
where $\Y \rightarrow S$ is flat with $d$-dimensional fibers.
For $i \ge 1$, let $s_i \in S(\C)$ and let $Y_i \subset V$ be
the fiber of $\Y \rightarrow S$ above $s_i$.
Let $\mu_i$ be a measure supported on $Y_i(\C)$.
If the $\mu_i$ converge weakly to a measure $\mu$ on $V(\C)$,
then the support of $\mu$ is contained in a $d$-dimensional Zariski closed
subvariety of $V$.
\end{lemma}

\begin{proof}
We may assume that $V=\PP^n$ and that $\Y \rightarrow \PP^n_S$ is
the universal family over a Hilbert scheme $S$.
Since $S$ is projective over $\C$,
we may pass to a subsequence to assume that the $s_i$
converge in the complex topology to $s \in S(\C)$.
By compactness of $\Y(\C)$, $\mu$ must be supported
on the fiber $\Y_s \subset V$.
\end{proof}

\begin{rem}
The hypotheses can be weakened.
It is enough to assume that the $Y_i$
form a ``limited family'' of closed subvarieties of $V$ of dimension $\le d$,
because then there are finitely many possibilities for their Hilbert
polynomials~\cite[Th\'eor\`eme~2.1]{fondements}.
The limited family condition holds, for instance,
if the $Y_i$ are reduced and equidimensional of dimension $d$,
and $\deg Y_i$
(with respect to some fixed embedding $V \hookrightarrow \PP^n$)
is bounded~\cite[Lemme~2.4]{fondements}.
\end{rem}

\begin{proof}[Proof of Theorem]
We may assume that $A = A_0 \times \G_m^n$.
Let $G \subset A$ be the connected component of the
group of translations preserving $X$.
We may assume $\dim G=0$, since otherwise we consider $X/G \subset A/G$
and use the easy fact that if $z_i \in A(\kbar)$ and $h(z_i) \rightarrow 0$
then the same is true for the images of the $z_i$ in $A/G$
(for any similarly defined height on $A/G$,
which is almost split).
We may enlarge $k$ so that $\Gamma \subset A(k)$.

If the theorem fails,
there exists a sequence $x_i=\gamma_i + z_i$ in $X(\kbar)$ converging
to the generic point of $X$,
with $\gamma_i \in \Gamma$ and $z_i \in A(\kbar)$
such that $h(z_i) \rightarrow 0$.
Passing to a subsequence and enlarging $k$,
we may assume that the closed points of $A$ corresponding to the $z_i$
are contained in a geometrically integral
subvariety $Z$ of $A$ defined over $k$,
and converge to the generic point of $Z$.
Let $Z_\sigma$ be the base extension of $Z$
by a fixed embedding $\sigma: \kbar \hookrightarrow \C$.
Let $O(z_i)$ denote the image of the $\Gal(\kbar/k)$-orbit of $z_i$
in $Z_\sigma(\C)$.
The Bogomolov conjecture and the equidistribution theorem
for almost split semiabelian varieties~\cite{chambertloir2}
imply that $Z$ is
the translate of a sub-semiabelian variety $B$
by a torsion point,
and that the uniform measures $\mu_i$ on $O(z_i)$
converge weakly as $i \rightarrow \infty$
to the measure $\mu_Z$ on $Z_\sigma(\C)$
that is the translate of Haar measure on the
maximal compact subgroup $B^1$ of $B_\sigma(\C)$.

If $Z$ is a point, we are done by Mordell-Lang.
Therefore assume $\dim Z>0$.
Let $\pi: A \rightarrow A/B$ be the projection.
If $\tau \in \Gal(\kbar/k)$, then
$\tau x_i - x_i = \tau z_i - z_i \in Z_\kbar - Z_\kbar = B_\kbar$,
so $\pi(\tau x_i)=\pi(x_i)$.

Applying an automorphism on the $\G_m^n$ side, we may assume that
$B=B_0 \times \G_m^r \times \{1\}^{n-r} \subset A_0 \times \G_m^n=A$,
where $B_0$ is a sub-abelian variety of $A_0$.
We embed $A$ in $\Abar:=A_0 \times (\PP^1)^n$ in the obvious way,
so that $\pi$ extends to a projective morphism
$\pibar: \Abar \rightarrow A_0/B_0 \times (\PP^1)^{n-r}$.
Let $\Xbar$ be the closure of $X$ in $\Abar$.
Restrict $\pibar$ to $\Xbar$ to obtain a projective morphism
$\pi_X: \Xbar \rightarrow \pibar(\Xbar)$
to the scheme-theoretic image.
Choose a dense open subset $U \subset \pibar(\Xbar)$
with $\pi_X^{-1}(U) \rightarrow U$ flat.
Discarding finitely many $x_i$,
we may assume that $\pi_X(x_i) \in U(\kbar)$ for all $i$.
Since $\dim G=0<\dim B$, the translation by elements in $B$
cannot all preserve $\Xbar$,
nor can they preserve the dense open subset $\pi_X^{-1}(U)$.
Hence the relative dimension $d$ of $\pi_X^{-1}(U) \rightarrow U$
satisfies $d < \dim B = \dim Z$.

Let $\Y \rightarrow S$ be the base extension
$\pi_X^{-1}(U) \times A \rightarrow U \times A$ of $\pi_X$.
Again $\Y$ is flat and projective over $S$ of relative dimension $d$.
We consider $\Y$ as a family of subvarieties of $\Abar$
via the $S$-morphism $\Y \rightarrow \Abar \times S$
where $\Y \rightarrow \Abar$ is the composition
$\Y = \pi_X^{-1}(U) \times A \hookrightarrow \Abar \times A \rightarrow \Abar$,
with the last map extending subtraction $A \times A \rightarrow A$.
The morphism
$\pi_X^{-1}(U) = \Xbar \times_{\pibar(\Xbar)} U \rightarrow \Xbar \times U$
is a closed immersion (since $\pibar(\Xbar)$ is separated over $k$).
Multiplying by $A$ and composing with
the product of $\Xbar \hookrightarrow \Abar$ by $U \times A$,
and the automorphism of
$\Abar \times U \times A$
given by ``subtracting'' the third coordinate from the first,
we find that
	$$\Y = \pi_X^{-1}(U) \times A \rightarrow \Xbar \times U \times A
		\rightarrow \Abar \times U \times A
		\rightarrow \Abar \times U \times A = \Abar \times S$$
is a closed immersion.

Let $s_i=\sigma(\pi_X(x_i),\gamma_i) \in S_\sigma(\C)$.
Let $Y_i$ be the fiber of $\Y$ above $s_i$,
considered as subvariety of $\Abar_\sigma$.
(Informally, $Y_i$ is $\pi_X^{-1}(\pi_X(x_i)) - \gamma_i$
base extended to $\C$.)
Using $\pi(\tau x_i)=\pi(x_i)$ for $\tau \in \Gal(\kbar/k)$,
we find that $\sigma(\tau z_i) \in Y_i(\C)$.
Thus $\mu_i$ is supported on $Y_i(\C)$.
The lemma applied to $\Y_\sigma \rightarrow S_\sigma$ with $V=\Abar_\sigma$
implies that the limit $\mu_Z$
is supported on a $d$-dimensional subvariety of $\Abar_\sigma$.

The Zariski closure of $B^1$ is an algebraic subgroup of $B_\sigma$
containing all its torsion, so $B^1$ is Zariski dense in $B_\sigma$.
Hence the support of $\mu_Z$ (which is a translate of $B^1$)
cannot be contained in a subvariety of dimension $d<\dim Z$.
This is a contradiction.
\end{proof}

\begin{rem} (Failure of equidistribution.)
Suppose $A$ is abelian,
$x_i \in X(\kbar) \cap \Gamma_{\epsilon_i}$
converge to the generic point of $X$,
and $\epsilon_i \rightarrow 0^+$.
The theorem tells us that $X$ must be a translated sub-abelian variety,
but it is {\em not} necessarily true that the uniform measures
on $O(x_i)$ (defined as above)
converge weakly to Haar measure on $X_\sigma(\C)$:
it could be that $x_i=n_i\gamma$ with $n_i \in \Z$, $\gamma \in \Gamma$
is a sequence of $k$-rational points converging in the complex topology,
for instance.
\end{rem}

\section{Conjecture}

In this section we formulate a conjecture which still includes
the generalized Bogomolov conjecture,
but now contains also the full Mordell-Lang conjecture for division points
on semiabelian varieties over number fields~\cite{mcquillan}.
The following approach, suggested to the author by E.~Hrushovski,
circumvents the problems with defining a canonical height
for semiabelian varieties.

Let $U$ be a geometrically integral quasi-projective variety
over a number field $k$, equipped with a morphism $f: U \rightarrow U$.
For integers $r \ge 1$,
let $f^r: U \rightarrow U$ denote the $r$-th iterate of $f$.
We assume the following condition on $(U,f)$:
\begin{itemize}
\item[$(\ast)$]	There exist a (logarithmic) Weil height
		$h: U(\kbar) \rightarrow \R$
	associated to some embedding $U \hookrightarrow \PP^n$,
	an integer $r \ge 1$, and real numbers $M>0$ and $c>1$
	such that $h(z)>M$ implies $h(f^r(z))>c h(z)$.
\end{itemize}
If $z \in U(\kbar)$, let $N(z)$ be the smallest integer $N \ge 1$
such that $h(f^N(z)) > M$, or $\infty$ if no such $N$ exists.
Northcott's finiteness theorem about the number of points
of bounded height and degree implies that $N(z)=\infty$
if and only if $z$ is {\em preperiodic} for $f$
(i.e., has finite orbit under the iterates of $f$).
For $i \ge 1$, let $z_i \in U(\kbar)$.
We say that $\{z_i\}_{i \ge 1}$ is a {\em sequence of small points}
if $N(z_i) \rightarrow \infty$ in $\PP^1(\R)$.

\begin{prop}
Suppose $(U,f)$ satisfies $(\ast)$,
with $h$, $r$, $M$, and $c$.
\begin{itemize}
	\item[1)] If $h'$ is another Weil height,
corresponding to another projective embedding,
then there exist $r'$, $M'$, and $c'$ as in $(\ast)$ for $h'$.
	\item[2)] For any such choices,
the notion of sequence of small points obtained is the same.
	\item[3)] If $(U,g)$ also satisfies $(\ast)$, and $fg=gf$,
then the notion obtained using $g$ is the same as that using $f$.
	\item[4)] If $(U',f')$ also satisfies $(\ast)$,
and $\psi: U \rightarrow U'$ satisfies $\psi f = f' \psi$,
then $\psi$ maps sequences of small points to sequences of small points.
\end{itemize}
\end{prop}

\begin{proof}
Without loss of generality, we may add a constant to the Weil heights
to assume that they always exceed~1.
(This allows us to drop a few constants in what follows.)

1)
Writing explicitly the isomorphisms
between the two embedded copies of $U$
shows that $h'(z) < e h(z)$
and $h(z) < e' h'(z)$ hold on $U(\kbar)$ for some $e,e'>1$.
Then for $M':=eM$, $c'=2$,
and for a sufficiently large multiple $r'=mr$ of $r$,
	$$h'(z)>M' \implies h(z)>M \text{ and }
		h'(f^{r'}(z)) > \frac{1}{e'} h(f^{r'}(z))
		> \frac{c^m}{e'} h(z)
		> \frac{c^m}{ee'} h'(z)
		> c' h'(z).$$

2)  Let $N'(z)$ be the $N$-function for a second choice $h'$, $M'$.
There exists an integer $p \ge 1$ independent of $z$,
such that
	$$h(z)>M \implies h'(f^{pr}(z)) > \frac{1}{e'} h(f^{pr}(z))
	> \frac{c^p}{e'} h(z) > \frac{c^p}{e'} M > M'.$$
Then $N'(z) \le N(z)+pr$ for all $z$.
The same argument bounds $N(z)$ by $N'(z)$ plus a constant.

3)  It suffices to show that $N_g$ (the $N$-function for $g$)
is bounded by a linear function of $N_f$.
(Then interchange $f$ and $g$.)
By~2), and by replacing each $r$ by a multiple,
we may assume that the $h$, $r$, $M$, and $c$
being used for $g$ are the same as for $f$.
Choose $d>1$ so that $h(f(z)) < d h(z)$ for all $z \in U(\kbar)$.
Let $\eta=N_f(z)$.
We may assume $\eta<\infty$.
For integers $m \ge (\log d/\log c) \eta$,
	$$h(g^{mr}(z)) > d^{-\eta} h(f^\eta(g^{mr}(z)))
		= d^{-\eta} h(g^{mr}(f^\eta(z)))
		> d^{-\eta} c^m h(f^\eta(z))
		\ge h(f^\eta(z))
		> M.$$
Hence $N_g(z) \le \lceil \log d / \log c \rceil r \cdot N_f(z)$.

4)
Choose $h'$ for $U'$.
Choose $\alpha>1$ so that $h'(\psi(z)) < \alpha h(z)$ holds
for all $z \in U(\kbar)$.
We may choose $M'$ for $U'$ so that $M'>\alpha M$.
Then $N'(\psi(z)) \ge N(z)$ for all $z \in U(\kbar)$.
\end{proof}

\begin{rem}
Condition $(\ast)$ is satisfied for $(U,f)$ if
there exists an integral projective variety $V$
containing $U$ as an open dense subset,
and an ample line bundle $\LL$ on $V$
such that $f$ extends to a morphism $\fbar: V \rightarrow V$ and
$\fbar^\ast \LL \tensor \LL^{\tensor -q}$ in $(\Pic V) \tensor \Q$
is {\em effective} for some $1< q \in \Q$.
(We say that $\MM \in (\Pic V) \tensor \Q$ is effective
if $\MM^{\tensor n}$ is represented by an effective Cartier divisor
for some $n \ge 1$.)
\end{rem}

\begin{rem}
Our situation is only slightly more general than that considered
in~\cite{callsilverman} and the introduction of~\cite{zhangsmall},
which consider projective varieties $X$ equipped with $f: X \rightarrow X$
and $\LL$ in $\Pic V$ or $(\Pic V) \tensor \R$
satisfying $f^\ast \LL = \LL^{\tensor d}$ for some $d>1$.
With our weaker assumptions,
we have apparently lost the ability
to define a reasonable ``canonical height'' using $f$,
but we have gained the ability to handle semiabelian varieties,
as we explain next.
\end{rem}

Let $A$ be a semiabelian variety over a number field $k$,
i.e., an algebraic group fitting into an exact sequence
	$$0 \rightarrow T \rightarrow A \rightarrow A_0 \rightarrow 0.$$
where $T_\kbar \isom \G_{m,\kbar}^r$ for some $r \ge 0$
and $A_0$ is an abelian variety.
We enlarge $k$ if necessary to assume that $T \isom \G_m^r$ over $k$.

For $m \in \Z$, let $[m]:A \rightarrow A$ be multiplication by $m$.
There is a compactification $\Abar$ to which $[m]$ extends,
equipped with effective line bundles $\LL_0$, $\LL_1$
with $\LL:=\LL_0 \tensor \LL_1$ ample,
such that $[m]^\ast \LL = \LL_0^{\tensor m^2} \tensor \LL_1^{\tensor m}$.
(See Section~1.1 of~\cite{mcquillan}, for example.)
The last remark above shows
that $(A,[m])$ satisfies $(\ast)$ for $m \ge 2$,
and the proposition shows that the resulting notion of
a sequence of small points depends only on $A$.

\begin{rem}
If $A$ is almost split, and $h$ is as in our introduction,
then $\{z_i\}_{i \ge 1}$ is a sequence of small points
if and only if $h(z_i) \rightarrow 0$.
\end{rem}

Let $\Gamma$ be a finitely generated subgroup of $A(\kbar)$,
and define the division group
	$$\Gamma' := \{\, x \in A(\kbar) \mid
	\text{ there exists $n \ge 1$ such that } nx \in \Gamma \,\}.$$
Finally, let $X$ be a geometrically integral closed subvariety of $A$.
In light of the many special cases that have been proven,
it seems reasonable to conjecture the following:

\begin{conj}
\label{mainconj}
For $i \ge 1$, suppose $x_i=\gamma_i+z_i \in X(\kbar)$
where $\gamma_i \in \Gamma'$
and $\{z_i\}_{i \ge 1}$ is a sequence of small points in $A(\kbar)$.
If $X$ is not a translate of a sub-semiabelian variety of $A$
by an element of $\Gamma'$,
then the $x_i$ are not Zariski dense in $X_\kbar$.
\end{conj}

\noindent
Equivalently, we could let
$B_\epsilon:=\{\,z \in A(\kbar) \mid N(z)>1/\epsilon \,\}$
and conjecture that for some $\epsilon>0$,
$X(\kbar) \cap (\Gamma' + B_\epsilon)$
is contained in a finite union $\bigcup Z_j$
where each $Z_j$ is a translate of a sub-semiabelian variety of $A_\kbar$
by an element of $\Gamma'$,
and $Z_j \subseteq X_\kbar$.
(Here $N$ is the $N$-function as above for $(A,[m])$,
for some Weil height $h$ and some $M>0$.)

\begin{rems}
The case where $\Gamma=0$ is a generalization
of the Bogomolov conjecture to semiabelian varieties.
Some foundations for a theory of canonical heights on semiabelian varieties
can be found in~\cite{chambertloir}.

The case where $A$ is almost split
and all the $\gamma_i$ are in $\Gamma$
is a restatement of our main theorem.
If we had a generalization of the equidistribution theorem
and the Bogomolov conjecture to general semiabelian varieties,
we could prove the case where $A$ is semiabelian and $\gamma_i \in \Gamma$
using our method.
\end{rems}

\section*{Acknowledgements}

I thank Ahmed Abb\`es and Jean-Beno\^{\i}t Bost
for many suggestions and corrections,
and the former for explaining the proof of the Bogomolov conjecture to me.
I also thank Ehud Hrushovski for the idea that one could
define the notion of a sequence of small points using Weil heights
even in situations where one does not have a good canonical height.
Finally I thank Shou-Wu Zhang for some interesting conversations on
related questions,
and Antoine Chambert-Loir for discussing his new results with me.

\end{document}